\documentclass[12pt]{article}
\usepackage[psamsfonts]{amssymb}

\advance\baselineskip by 22.4pt plus1pt
\def\bra#1{\left\langle #1\right|}
\def\ket#1{\left| #1\right\rangle}
\def\hs#1#2{\left\langle #1|#2\right\rangle}

\def\Mat{{\rm Mat}\,}

\def\nn{\nonumber}
\newcommand{\IZ}{{\mathbb Z}}
\newtheorem{proposition}{Proposition}[section]  
\newcommand{\bprop}{\medskip\begin{proposition} ~~\\ \it}
\newcommand{\eprop}{\end{proposition} \hfill \medskip \\}
\def\b#1{{\mathbb #1}}
\def\bra#1{\langle #1 \vert}
\def\c#1{{\cal #1}}

\def\Dirac{{\raise0.09em\hbox{/}}\kern-0.69em D}

\def\ket#1{\vert #1 \rangle}

\def\exterior{{{\raise0.2em\hbox{$\scriptstyle\bigwedge$}}{}}}
\def\k{\kern-.1em\mathbin{,}\kern-.1em}
\def\hk{\kern.12em\raise-1em\hbox{$\hat{\raise1em\hbox{,}}$}\kern.12em}
\newcommand{\initiate}{\setcounter{equation}{0}}
\begin{document}

\title{Twisted Configurations \\ over \\
 Quantum Euclidean Spheres}

\author{Giovanni Landi$\strut^{1},$\ John Madore$\strut^{2,3}$,\\[15pt]
             $\strut^{1}$Dipartimento di Scienze Matematiche, Universit\`a
             di Trieste\\
             Via Valerio 12/b, I-34127 Trieste \\[5pt]
             $\strut^{2}$Laboratoire de Physique Th\'eorique\\
             Universit\'e de Paris-Sud, B\^atiment 211, F-91405 Orsay
\and    $\strut^{3}$Max-Planck-Institut f\"ur Physik\\
             F\"ohringer Ring 6, D-80805 M\"unchen}

\date{}

\maketitle

\begin{abstract}
     	We show that the relations which define the algebras of 
	the quantum Euclidean planes $\b{R}^N_q$ can be expressed in terms of
        projections provided that the unique central element, the radial
        distance from the origin, is fixed.  The resulting reduced algebras
        without center are the quantum Euclidean spheres $S^{N-1}_q$.
	The projections $e=e^2=e^*$ are elements in 
	$\Mat_{2^n}(S^{N-1}_q)$, with $N=2n+1$ or $N=2n$, and can be regarded 
	as defining modules of sections of $q$-generalizations of monopoles,
	instantons or more general twisted bundles over the spheres.
	We also give the algebraic definition of normal and cotangent
	bundles over the spheres in terms of canonically defined
     	projections in $\Mat_{N}(S^{N-1}_q)$.
\end{abstract}

\vfill
May 2002 \\
Revised: June 2002

\newpage

\initiate
\section{Introduction}

Let $\c{A}$ be an associative algebra defined formally in terms of a
set of $N$ generators $x^i$ and a set of relations $R(x^i)$. Suppose
that for a fixed value of each of the parameters which enter in the
definition of the $R(x^i)$ some of the relations reduce to the
constraint that the algebra be commutative. We refer to these as
commutation relations. The remaining relations become constraints in
the commutative limit and define the noncommutative generalization of
a submanifold of $\b{R}^N$. We suppose for simplicity that there is
only one relation of this form.  The commutative limit will describe then a
submanifold $V$ of dimension $N-1$ embedded in $\b{R}^N$. If one
introduce the moving frame $dx^i$ on $\b{R}^N$ then the latter
acquires the structure of a flat differential manifold. This defines
by the embedding, that is, by the relations, a moving frame
$\theta^\alpha$ on $V$, which one can always choose so that locally
the last component is normal to $V$. When $V$ is parallelizable this
choice can be made globally and the module of sections $\c{T}$ of the
cotangent bundle $T^*(V)$ is free. In all cases the embedding defines
a splitting of the module of sections $\c{A}^N$ of $T^*(\b{R}^N)$ into
a direct sum
\begin{equation}
\c{A}^N = \c{T} \oplus \c{N}
\end{equation}
of $\c{T}$ and a free module $\c{N}$ of rank one.  The metric and
the frame on $V$ are determined by the embedding relations.  The
construction is most elegantly described in terms of projections in the matrix
algebra $\Mat_{N}(\c{A})$, i.e. elements $e\in \Mat_{N}(\c{A})$ such
that $e^2=e=e^*$.  It is a consequence of the Serre-Swan theorem (cf.
\cite{Con94,Lan97,Mad99,Gra00}) that to every
vector bundle over $V$ corresponds an equivalence class of
projections in a matrix algebra
$\Mat_{r}(\c{A})$, for a suitable $r$. The module of sections of the
bundle is a projective module of finite type over the algebra $\c{A}$ and
conversely any such a module can be realized as the module of sections of a
bundle; and these modules are naturally characterized by projections. This
correspondence quite naturally generalizes the notion of a vector bundle over a
noncommutative algebra.

We shall see that in the case of
the quantum Euclidean spheres $S^{N-1}_q$ a projection
$e\in\Mat_{2^n}(S^{N-1}_q)$,  with $N=2n$ or $N=2n+1$, can be so
chosen  that the
relations which define it ($e^2=e=e^*$) are equivalent to the relations
$R(x^i)$ which define the algebra $\c{A}$ of the spheres.
As already mentioned, the projections can be regarded as defining modules
of  sections of bundles over the spheres which then will be
$q$-generalizations  of monopoles, instantons or more general twisted
configurations.  
As in the classical situation, the bundles
should be characterized by integer valued topological charges (Chern numbers);
work on these is in progress and will be reported somewhere else.

In Section~\ref{se:frt} some basic 
results~\cite{FadResTak89} about the
general $N$-dimensional quantum Euclidean spaces $\b{R}^N_q$ and quantum 
Euclidean spheres $S^{N-1}_q$ are reviewed.  In Section~\ref{se:ins} we shall 
introduce the projections which determine (and are determined by) the spheres 
relations and which, in turn, define twisted configurations over the spheres
$S^{N-1}_q$. 
We shall describe in detail the cases $N=3,4,5,6$, that is,
the spheres of  dimensions $2,3,4$ and $5$ respectively and
outline the general case. Not  surprisingly, `even and odd
dimensional' spheres will behave differently.  In
Section~\ref{se:nct} we construct the projections  for the
normal and cotangent bundles over the spheres $S^{N-1}_q$. 
The final remarks concern some preliminary results on the computation of
topological charges.

\initiate

\section{Quantum Euclidean Planes and Spheres}
\label{se:frt}

In this section some basic results~\cite{FadResTak89} (see also 
\cite{CerFioMad00a}) about the general $N$-dimensional quantum Euclidean 
spaces $\b{C}^N_q$ and $\b{R}^N_q$ and spheres $S^{N-1}_q$ are reviewed.  
We start with the matrix $\hat{R}$ for the quantum group
$SO_q(N,\b{C})$.  It is a symmetric $N^2 \times N^2$ matrix and its
main property
is that it satisfies the braid relation.  It admits a projector decomposition:
\begin{equation}
\hat R = q P_{(s)} - q^{-1} P_{(a)} + q^{1-N} P_{(t)}.   \label{projectorR}
\end{equation}
where the $P_{(s)}$, $P_{(a)}$, $P_{(t)}$ are $SO_q(N)$-covariant
$q$-deformations of the symmetric trace-free, antisymmetric and trace
projectors respectively.  The projector $P_{(t)}$ projects onto a
one-dimensional sub-space and can be written in the form
$P_{(t)}^{ij}{}_{kl} = (g^{mn}g_{mn})^{-1} g^{ij}g_{kl}$. This leads
to the definition of a deformed metric matrix, the $N \times N$ matrix
$g_{ij}$ given by 
\begin{equation}
g_{ij} = q^{-\rho_i} \delta_{i,-j}.
\end{equation}
which is a $SO_q(N)$-isotropic tensor. If $n$ is the rank of $SO(N,\b{C})$, i.e. the
integer part in $ N / 2$, the indices take the values
$i=-n,\ldots,-1,0,1,\ldots, n$ for $N=2n+1$, and $i=-n,\ldots,-1,
1,\ldots, n$ for $N=2n$.  Moreover, we have introduced the notation
\begin{equation}
\rho_i=\left\{\begin{array}{ll}
(n-\frac{1}{2},\ldots,\frac{1}{2},0,-\frac{1}{2},\ldots,\frac{1}{2}-n)
& \hbox{ for $N = 2n +1 $,}\\
(n-1,\ldots,0,0,\ldots,1-n) & \hbox{ for $N = 2n$.}
\end{array}\right.
\end{equation}
The metric and the braid matrix satisfy the `$gTT$'
relations~\cite{FadResTak89}
\begin{equation}
g_{il}\,\hat R^{\pm 1}{}^{lh}_{jk} =
\hat R^{\mp 1}{}^{hl}{}_{ij}\,g_{lk}, \qquad
g^{il}\,\hat R^{\pm 1}{}{}_{lh}^{jk} =
\hat R^{\mp 1}{}{}_{hl}^{ij}\,g^{lk}.                          \label{gRrel}
\end{equation}

With the help of the projector $P_{(a)}$, the $N$-dimensional quantum
Euclidean space is defined as the associative algebra $\b{C}^N_q$
generated by elements $\{x^i\}_{i=-n,\ldots,n}$ with relations
\begin{equation}
P_{(a)}{}^{ij}{}_{kl} x^k x^l=0.                              \label{xrel}
\end{equation}
or, more explicitly~\cite{Ogi92a}
\begin{equation}
x^i x^j = q x^j x^i \hbox{ for } i<j, i \neq -j, \qquad
[x^i, x^{-i}]=\left\{\begin{array}{ll}
k \omega_{i-1}^{-1} r^2_{i-1} &\hbox{ for } i > 1\\
0 & \hbox{ for $i=1$, $N=2n$},\\
h r_0^2 & \hbox{ for $i=1$, $N=2n+1$} .
\end{array}\right.                                         \label{explicitx}
\end{equation}
We use the notation
$\omega_i=q^{\rho_i}+q^{-\rho_i}$,
$h=q^{\frac{1}{2}}-q^{-\frac{1}{2}}$,
$k=q-q^{-1}$ and
\begin{equation}
r^2_i=\sum_{k,l=-i}^i g_{kl} x^k x^l, \qquad
i \ge 0 \mbox{ for $N=2n+1$},
\:  i \ge 1 \mbox{ for $N=n$} .                              
\label{defr}
\end{equation}
The last element $r^2 \equiv r^{2}_{n}$ can be shown to be central.

For $q \in \b{R}^+$ a conjugation $(x^i)^*= x^j g_{ji}$ can be
defined on $\b{C}^N_q$ to obtain what is known as the quantum real
Euclidean space $\b{R}^N_q$. The relations (\ref{explicitx}) can be
used to derive analogous ones
for the variables $\{(x^i)^* , x^j\}$,
\begin{eqnarray}
&& x^i x^j = q x^j x^i \hbox{ for } i<j,  \qquad
x_i^* x_j = q x_j x_i^*, ~~~ i\neq j
\nonumber \\
&& ~ \nonumber \\
&& [x^i, (x^i)^*]=\left\{\begin{array}{ll}
[ (1 - q^{-2}) /  (1 + q^{-2\rho_{i-1}} )]
 ~ r^2_{i-1} &\hbox{ for } i > 1\\
0 & \hbox{ for $i=1$, $N=2n$},\\
(1-q^{-1}) ~ r_0^2 & \hbox{ for $i=1$, $N=2n+1$} .
\end{array}\right.                                         
\label{explicitxbis}
\end{eqnarray}
The central element $r^2$ can be
written as
\begin{equation}
r^2=\sum_{k=-n}^n (x^k)^* x^k = q^{-2\rho_n} x^n (x^n)^* + \cdots + 
(x^n)^* x^n,                              \label{radius}
\end{equation}
with similar expressions for the elements $r^2_i$, $i \ge 0 \mbox{ 
for 
$N=2n+1$}, \:  i \ge 1 \mbox{ for $N=n$}$,
\begin{equation}
r^2_i=\sum_{k=-i}^i (x^k)^* x^k = q^{-2\rho_i} x^i (x^i)^* + 
\cdots + (x^i)^* x^i ,
\label{defrbis}
\end{equation}
By fixing the value of $r^2$ we get the quantum Euclidean sphere 
$S^{N-1}_q$
of the corresponding radius.  Thus, the quantum Euclidean sphere 
$S^{N-1}_q$ is naturally considered as a quantum subspace of the quantum
Euclidean space $\b{R}^N_q$ and the algebra of functions on the sphere is a quotient
of the algebra of functions on the quantum Euclidean space, the quotient being taken by
the ideal generated by the relation that fixes the radius. Furthermore, 
the spheres $S^{N-1}_q$ are invariant
under the action of the quantum groups $SO_q(N)$ on them. 

It is easy to see that the spheres $S^{N-1}_q$ have a $S^1$ 
worth of classical points. Indeed, let $r\in\b{R}\setminus \{0\}$ be the 
radius of the sphere.
Then, with $\lambda \in \b{C}$ such that $|\lambda|^2=1$, there is a 
family 
of $1$-dimensional representations (characters) of the algebra
$S^{N-1}_q$ given by
\begin{eqnarray}\label{1reps}
&& \tau_{\lambda}(1) = 1, 
~~\tau_{\lambda}(x^{n}) = {r \over \sqrt{1+q^{-2\rho_n}}} ~\lambda ,
~~\tau_{\lambda}((x^{n})^*) = {r \over \sqrt{1+q^{-2\rho_n}}} ~\bar{\lambda}~, \\ 
\nonumber
&& \tau_{\lambda}(x^{i}) = \tau_{\lambda}((x^{i})^*) = 0 ~, 
\end{eqnarray}
for $i=0,1,\ldots, n-1$ or $i=1,\ldots, n-1$ according to whether $N=2n+1$ 
or 
$N=2n$. Clearly these representations will yield traces on the algebras $S^{N-1}_q$.

\initiate
\section{Twisted configurations over $S^{N-1}_q$}\label{se:ins}

We shall now introduce a hermitian idempotent or projection, an
element $e \in \Mat_{2^n}(S^{N-1}_q)$ satisfying the conditions 
$e = e^2 = e^*$ which determine the 
structure of the algebra
$S^{N-1}_q$. The star here is the formal adjoint in the complete
matrix algebra and $N$ is even or odd. These projections, in turn, can
be thought of as defining modules of sections of bundles over the 
spheres $S^{N-1}_q$ which
are
$q$-generalizations of monopoles and instanton bundles or more general twisted
configurations. We shall use a matrix trace together with the trace determined by the
representations (\ref{1reps}) to compute the rank which, for the projector in
question,  will turn out to be equal to $2^{n-1}$. 
At the moment we are unable to compute the topological charges in general; work on this
is in progress and will be reported elsewhere.  We shall describe in detail the cases
$N=3,4,5,6$, that is the spheres of dimensions
$2,3,4$ and $5$ respectively. We shall also outline the general case
which in principle can be obtained using the same techniques.
\subsection{Monopoles on the Euclidean sphere $S^2_q$}

With a suitable rescaling of the 
generators the Euclidean sphere $S^2_q$ can be identified 
with the so called equator sphere of Podle{\'s} \cite{Pod87}. 
For its presentation we have the following generators
\begin{equation}\label{gen2}
x^i = (\sqrt q x_1^*, x_0, x_1) ,
\end{equation}
with $(x_0)^* = x_0$,
and commutation relations (\ref{explicitxbis}) given by
\begin{eqnarray}\label{cr2}
&& x_0 x_1 = q x_1 x_0, \qquad x_1^* x_0 = q x_0 x_1^*, \nonumber \\
&& [x_1, x_1^*] = (1-q^{-1}) x_0^2 .
\end{eqnarray}
These commutation relations give for the central element $r^2$ the
equivalent expressions
\begin{eqnarray}\label{rad2}
r^2 &=& q x_1 x_1^* + x_0^2 + x_1^* x_1 \nonumber \\
&=& q x_0^2 + (1+q) x_1^* x_1 = q^{-1} x_0^2 + (1+q) x_1 x_1^* .
\end{eqnarray}
Then, a straightforward computation yields that the element $e\in\Mat_{2}(S^2_q)$ given
by
\begin{equation}\label{projs2}
e = \frac 12 \left(\begin{array}{cc}
1+ q^{-1 / 2} ~r^{-1}~ x_0  ~&~ (1+q)^{1/2} ~r^{-1}~ x_1 \\[3pt]
(1+q)^{1/2} ~r^{-1}~ x^*_1 &1-q^{1/2} ~r^{-1}~ x_0 \end{array}\right)
\end{equation}
which is hermitian by construction, is also idempotent, $e^2=e$, if
and only if all the
relations (\ref{cr2}) and (\ref{rad2}) which define $S^2_q$ are satisfied.

When the  Euclidean sphere $S^2_q$ is identified with 
Podle{\'s} equator sphere the projector (\ref{projs2}) 
coincides with the one found in \cite{BrzMaj00} where projectors for all Podle{\'s}
spheres were constructed (a projector on the so called Podle{\'s} standard sphere had
already been  constructed in \cite{HajMaj98}).
\\ The projection (\ref{projs2}) (and the corresponding vector bundle over $S^2_q$)
is of rank $1$,
\begin{equation}
{\rm rank}(e) =: \tau_\lambda \circ Tr(e) =  \tau_\lambda 
\Big(1 + { 1 \over 2 r} ~q^{-1/2}~ (1-q) ~x_0 \Big) = 1~,
\end{equation}
where we have used the $1$-dimensional representations (\ref{1reps}) and $Tr$
denotes a matrix trace. In order to compute the topological charge of the bundle 
we need a  cyclic
$0$-cocycle, i.e. a trace $\tau^1$, on the reduced algebra $\bar{S^2_q} = S^2_q / 
\b{C} 1$ .
This computation is an example of the pairing between K-theory and cyclic cohomology
\cite{Con94}. In the present situation the cyclic cocycle one needs  is of degree zero since
for $S^2_q$ all cyclic cohomology is in $HC^0(S^2_q)$ while the even cyclic cohomology
$HC^n(S^2_q), n > 0$, is the image of the periodicity operator \cite{MasNakWat91}.
The trace $\tau^{1}$ was found in \cite{MasNakWat91} and on the generators
(\ref{gen2}) gives in particular 
\begin{equation}\label{2tra}
\tau^1(x_0)) = - 2 r ~q^{1/2} ~{ 1 \over 1 - q } ~.
\end{equation}
Then, with $[Tr(e)] \in \bar{S^2_q}$,
the topological charge of the projection
$e$ is found to be  
\begin{equation}\label{topcar}
\tau^{1} ([Tr(e)]) = \tau^{1} \Big( ~{ 1 \over 2 r} ~q^{-1/2}~ (1-q) ~x_0 \Big)
= -1~, 
\end{equation} 
as it should be given the analogous computations in 
\cite{BrzMaj98} and \cite{HajMaj98} for the computations of  
topological charges of
$q$-monopoles on Podle{\'s} spheres. The fact that the pairing in (\ref{topcar})
gives an integer number is a consequence of a noncommutative index theorem
\cite{Con94} since the trace $\tau^{1}$ is the character of a Fredholm module
\cite{MasNakWat91}. The trace (\ref{2tra}) is singular in the classical limit
$q=1$. In the latter limit  the bundle
corresponding to the projection (\ref{projs2}) is the
monopole bundle over the classical sphere $S^2$ of topological charge
(first Chern number)
equal to $-1$ which is computed by integrating a $2$-form on $S^2$ \cite{Lan00}
(classically the relevant cyclic cocycle is in degree two and is associated with the
volume form of the sphere).

\subsection{Instantons on the Euclidean sphere $S^4_q$}
In this section we generalize the projector (\ref{projs2}) on the sphere $S^2_q$ to a
projector on the sphere  $S^4_q$ and indicate how to further generalize it to any even
quantum Euclidean sphere. \\ 
The generators of the algebra  $S^4_q$ are the five elements
\begin{equation}
x^i = (q^{3/2} x_2^*, q^{1/2} x_1^*, x_0, x_1, x_2) ,
\end{equation}
with $(x_0)^* = x_0$.
Then the commutation relations (\ref{explicitxbis}) are
\begin{eqnarray}\label{cr4}
&& x_i x_j = q x_j x_i, ~~ i<j ; ~~~~~ x_i^* x_j = q x_j x_i^*, ~~~ 
i\neq j , \nonumber \\
&& [x_1, x_1^*] = (1-q^{-1}) x_0^2 , \nonumber \\
&& [x_2, x_2^*] = q^{-1}(1-q^{-1}) ~\Big( q x_1 x_1^* + x_0^2 + x_1^* x_1 \Big)
\nn \\
&& ~~~~~~~~~ = q^{-1}(1-q^{-1}) ~\Big( q x_0^2 + (1+q) x_1^* x_1 \Big) \nn \\
&& ~~~~~~~~~ = q^{-1}(1-q^{-1})~\Big( q^{-1} x_0^2 + (1+q) x_1 x_1^* \Big).
\end{eqnarray}
These commutation relations give for the central element $r^2$ the
equivalent expressions
\begin{eqnarray}
r^2 &=& q^3 x_2 x_2^* + q x_1 x_1^* + x_0^2 + x_1^* x_1 + x_2^* x_2
\nonumber \\
&=& (1+q^3) x_2^* x_2 + (1+q^3) x_1^* x_1 + q {1+q^3 \over 1+q}
x_0^2 \nonumber \\
&=& (1+q^3) x_2 x_2^* + q^{-2} (1+q^3) x_1 x_1^* + q^{-3} {1+q^3
\over 1+q} x_0^2 .
\end{eqnarray}
To alleviate the notation we now consider the case 
$r^2 = [ (1+q^3) / (1+q) ] \cdot 1$. Equivalently we could have
suitably rescaled the generators.  The previous sphere relations
reduce  to
\begin{eqnarray}\label{rad4}
1 &=& (1+q) x_2^* x_2 + (1+q) x_1^* x_1 + q x_0^2 \nonumber \\
&=& (1+q) x_2 x_2^* + q^{-2} [~ (1+q) x_1 x_1^* + q^{-1} x_0^2 ~].
\end{eqnarray}
Then, a straightforward computation yields that the 
element $e\in\Mat_{4}(S^4_q)$ given by
\begin{equation}\label{projs4}
e = \frac 12 \left(\begin{array}{cccc}
1 + q^{-3/2} x_0 & q^{-1} (1+q)^{1/2} x_1 & (1+q)^{1/2} x_2 & 0 \\[3pt]
q^{-1} (1+q)^{1/2} x^*_1 & 1 - q^{-1/2} x_0 & 0 & (1+q)^{1/2} x_2 \\[3pt]
(1+q)^{1/2} x^*_2 & 0 & 1 - q^{-1/2} x_0 & -(1+q)^{1/2} x_1 \\[3pt]
0 & (1+q)^{1/2} x^*_2 & -(1+q)^{1/2} x^*_1 & 1 + q^{1/2} x_0
\end{array}\right)
\end{equation}
which is hermitian by construction, is also idempotent, $e^2=e$, if
and only if all the relations (\ref{cr4}) and and (\ref{rad4}) which define
$S^4_q$ are satisfied.  
For a generic value of the radius $r$ the corresponding idempotent
can be easily guessed from the
previous expression and the analogue one (\ref{projs2}) for the
sphere $S^2_q$: one has simply to multiply the algebra generators in
(\ref{projs4}) by the constant $[ (1+q^3) / (1+q) ]^{1/2} r^{-1}$.

\noindent
It is worth stressing the way the relations (\ref{rad2}) for the
sphere $S^2_q$ are incorporated in
the ones (\ref{rad4}) for the sphere $S^4_q$. This fact helped in
constructing the
projection (\ref{projs4}) for $S^4_q$ by a suitable use of the projection
(\ref{projs2}) for $S^2_q$. Indeed, the projection $e_{(4)}$ in (\ref{projs4})
is related to the analogous one $e_{(2)}$ in (\ref{projs2}) by
\begin{equation}\label{projs4bis}
e_{(4)} = \frac 12 \left(\begin{array}{cc}
\b{I}_2 + q^{-1} u_{(2)} &  (1+q)^{1/2} x_2 ~\b{I}_2 \\[3pt]
(1+q)^{1/2} x^*_2 ~\b{I}_2 & \b{I}_2 - u_{(2)}
\end{array}\right) ,
\end{equation}
with the hermitian isometry $u_{(2)}$ given by $u_{(2)} = 2 e_{(2)} - 1$.
By applying this `inductive construction' one can
obtain the projections for the higher dimensional even spheres.

\noindent
By using the $1$-dimensional representations (\ref{1reps}) and a matrix trace $Tr$
one finds that the projection (\ref{projs4}) (and
the corresponding vector bundle over $S^4_q$) is of rank $2$,
\begin{equation}
{\rm rank}(e) =: \tau_\lambda \circ Tr(e) = \tau_\lambda \Big(2 + {1 \over 2}
~q^{-3/2} ~(1 - q )^2 ~x_0 \Big) = 2 ~.
\end{equation}
As it happens for the sphere $S^2_q$ described previously, in order to compute the
topological charge of the bundle, we need a  cyclic $0$-cocycle,
i.e. a trace $\tau^1$, on the reduced algebra $\bar{S^4_q} = S^4_q / 
\b{C} 1$, which
needs to be combined with the matrix trace. Now, the equivalence 
class of  $Tr(e)$ in
$\bar{S^4_q}$ is given by
\begin{equation}
[Tr(e)] = {1 \over  2} ~q^{-3/2} ~(1-q)^2 ~x_0 ~.
\end{equation}
and we expect that the trace $\tau^1$ should yield $-1$ on the previous 
expression. The reason to expect a value $-1$ for the topological charge is 
that in the classical limit
$q=1$, the bundle corresponding to the projection (\ref{projs4}) is the instanton
bundle over the classical sphere $S^4$ of topological charge (second Chern
number) equal to $-1$;  this is computed by integrating a $4$-form on $S^4$
\cite{Lan00}. As for the sphere $S^2_q$, the
trace $\tau^1$ on $S^4_q$ will be then singular in the classical limit. As for now we
have been unable to find such a singular trace on the sphere $S^4_q$.

\subsection{Solitons on odd quantum Euclidean spheres}
As already mentioned, odd quantum Euclidean spheres behave in a different manner
than the even ones in so that the corresponding projectors are trivial (i.e. they
correspond to trivial bundles). On $S^3_q$ this is a consequence of the fact that for the
K-theory group one has $K_0(S^3_q) = \IZ$
\cite{MasNakWat90} and for higher values of $N$ one expects a similar result.  

\subsubsection{The Euclidean sphere $S^3_q$} 
With a suitable rescaling of the generators  the
Euclidean sphere $S^3_q$ can be identified  with the `quantum sphere' $SU_q(2)$
of \cite{Wor87}.  The present generators of the algebra are four elements
\begin{equation}
x^i = (q x_2^*, x_1^*, x_1, x_2) 
\end{equation}
with commutation relations (\ref{explicitxbis}) 
\begin{eqnarray}\label{cr3}
&& x_i x_j = q x_j x_i, ~~ i<j ; ~~~~~ x_i^* x_j = q x_j x_i^*, ~~~ 
i\neq j , \nonumber \\
&& [x_1, x_1^*] = 0 , \qquad [x_2, x_2^*] = (1-q^{-2})x_1 x_1^* .
\end{eqnarray}
The central radial variable can be written as
\begin{eqnarray}
r^2 &=&q^2 x_2 x_2^* + x_1 x_1^* + x_1^* x_1 + x_2^* x_2
\nonumber \\
&=& (1+q^2) \Big( x_2^* x_2 + x_1^* x_1 \Big) = (1+q^2) \Big( x_2 x_2^* + 
q^{-2} x_1 x_1^*  \Big) .
\end{eqnarray}
To alleviate the notation we have chosen here the case 
$r^2 = (1+q^2) \cdot 1$ so that the previous sphere relations reduce to
\begin{equation}\label{rad3}
1 = x_2^* x_2 + x_1^* x_1 = x_2 x_2^* + q^{-2} x_1 x_1^* .
\end{equation}
Then,  
straightforward computations yield  that the
element $e\in\Mat_{4}(S^3_q)$ given by
\begin{equation}\label{projs3}
e = \frac 12 \left(\begin{array}{cccc}
1 & q^{-1} x_1 & x_2 & 0 \\[3pt]
q^{-1} x^*_1 & 1 & 0 & x_2 \\[3pt]
x^*_2 & 0 & 1 & -x_1 \\[3pt]
0 & x^*_2 & -x^*_1 & 1
\end{array}\right)
\end{equation}
which is hermitian by construction, is also idempotent, $e^2=e$, if
and only if all the
relations (\ref{cr3}) and and (\ref{rad3}) which define $S^3_q$ are satisfied.  
For a generic value of the radius $r$ the corresponding idempotent
can be easily guessed from the
previous expression: one has simply to multiply the algebra generators in
(\ref{projs3}) by the constant
$(1+q^2)^{1/2} r^{-1}$.

\noindent
The projection (\ref{projs3}) (and
the corresponding vector bundle over $S^3_q$) is of rank $2$,
\begin{equation}
{\rm rank}(e) = \tau_\lambda \circ  Tr(e) = \tau_\lambda (2) = 2~, 
\end{equation}
with the representation $\tau_\lambda$ given in (\ref{1reps}). By direct
computation, one checks that the Chern Character of the projection (\ref{projs3})
vanished identically.

\subsection{Solitons on the Euclidean sphere $S^5_q$}
The generators of the algebra $S^5_q$ are the six elements
\begin{equation}
x^i = (q^2 x_3^*, q x_2^*, x_1^*, x_1, x_2, x_3) .
\end{equation}
with commutation relations (\ref{explicitxbis}) given by
\begin{eqnarray}\label{cr5}
&& x_i x_j = q x_j x_i, ~~ i<j ; ~~~~~ x_i^* x_j = q x_j x_i^*, ~~~ 
i\neq j , \nonumber \\
&& [x_1, x_1^*] = 0 , \qquad [x_2, x_2^*] = (1-q^{-2})x_1 x_1^* , \nonumber \\
&& [x_3, x_3^*] = {1-q^{-2} \over 1+q^2} ~[q^2 x_2 x_2^* + x_1 x_1^* 
+ x_1^* x_1 +
x_2^* x_2], \nonumber \\
&& ~~~~~~~~~ = (1-q^{-2}) (x_2^* x_2 + x_1^* x_1) = (1-q^{-2})
(x_2 x_2^* + q^{-2} x_1 x_1^*) .
\end{eqnarray}
It follows that the radial element can be written as
\begin{eqnarray}
r^2 &=& q^4 x_3 x_3^* + q^2 x_2 x_2^* + x_1 x_1^* + x_1^* x_1 + x_2^* 
x_2 + x_3^* x_3 \nonumber \\ 
&=&
 (1+q^4) [ x_3^* x_3 + x_2^* x_2 + x_1^* x_1 ] = (1+q^4) [ x_3 
x_3^* + q^{-2} ( x_2
x_2^* + q^{-2} x_1 x_1^* )] .
\end{eqnarray}
For the moment we consider the case $r^2 = (1+q^4) \cdot 1$.
Then, the previous sphere relations
reduce to
\begin{equation}\label{rad5}
1 = x_3^* x_3 + x_2^* x_2 + x_1^* x_1
= x_3 x_3^* + q^{-2} ( x_2 x_2^* + q^{-2} x_1 x_1^* ) .
\end{equation}
And again, 
straightforward computations yield that the 
element $e\in\Mat_{8}(S^5_q)$ given by
\begin{equation}\label{projs5}
e = \frac 12 \left(\begin{array}{cccccccc}
1 & q^{-2}x_1 & q^{-1}x_2 & 0 & x_3 & 0 & 0 & 0 \\[3pt]
q^{-2}x^*_1 & 1 & 0 & q^{-1}x_2 &  0 &x_3 & 0 & 0 \\[3pt]
q^{-1}x^*_2 & 0 & 1 & -q^{-1}x_1 & 0 & 0 & x_3 & 0 \\[3pt]
0 & q^{-1}x^*_2 & -q^{-1}x^*_1 & 1 & 0 & 0 & 0 & x_3 \\[3pt]
x^*_3 & 0 & 0 & 0 & 1 & -q^{-1}x_1 & -x_2 & 0 \\[3pt]
0 & x^*_3 & 0 & 0 & -q^{-1}x^*_1 & 1 & 0 & -x_2\\[3pt]
0 & 0 & x^*_3 & 0 & -x^*_2 & 0 & 1 & x_1 \\[3pt]
0 & 0 & 0 & x^*_3 & 0 & -x^*_2 & x^*_1 & 1\\
\end{array}\right)
\end{equation}
which is hermitian by construction, is also idempotent, $e^2=e$, if
and only if all the
relations (\ref{cr5}) and and (\ref{rad5}) which define $S^5_q$ are satisfied.  
For a generic value of the radius $r$ the corresponding idempotent
can be easily guessed from the
previous expression: one has simply to multiply the algebra generators in
(\ref{projs5}) by the constant
$(1+q^4)^{1/2} r^{-1}$.

\noindent
Again, as for the even spheres, one should notice the way the
relations (\ref{rad3}) for the sphere $S^3_q$ are incorporated in
the ones (\ref{rad5}) for the sphere $S^5_q$. This fact helped in
constructing the
projection (\ref{projs5}) for $S^5_q$ by a suitable use of the projection
(\ref{projs3}) for $S^3_q$.
Indeed, the projection $e_{(5)}$ in (\ref{projs5})
is related to the analogous one $e_{(3)}$ in (\ref{projs3}) by
\begin{equation}\label{projs5bis}
e_{(5)} = \frac 12 \left(\begin{array}{cc}
\b{I}_4 + q^{-1} u_{(3)} & x_3 ~\b{I}_4 \\[3pt]
x^*_3 ~\b{I}_4 & \b{I}_4 - u_{(3)}
\end{array}\right) ,
\end{equation}
with the hermitian isometry $u_{(3)}$ given by $u_{(3)} = 2 e_{(3)} - 1$.
By applying this `inductive construction' one can
obtain the projections for the higher dimensional odd spheres.

\noindent
The projection (\ref{projs5}) (and
the corresponding vector bundle over $S^5_q$) is of rank $4$,
\begin{equation}
{\rm rank}(e) = \tau_\lambda \circ Tr(e) = \tau_\lambda (4) = 4.
\end{equation}
As for the $S^3_q$ and the corresponding projection (\ref{projs3}), the Chern
Character of the projection (\ref{projs5}) vanished identically as well and
the corresponding bundle is again trivial.

\initiate
\section{Normal and cotangent bundles}\label{se:nct}

We shall now introduce a projection in a suitable matrix algebra which
determines natural cotangent bundles over the spheres $S^{N-1}_q$. 
Consider then
the following `vector valued function' on $S^{N-1}_q$
\begin{equation}                                           \label{norbra}
\bra{\perp\!}=\left\{\begin{array}{ll}
{1 \over r} ~\Big( q^{-\rho_{n}} ~x_{n} , q^{-\rho_{n-1}} ~x_{n-1},
\ldots, q^{-\rho_{1}} ~x_{1} , x_{0} , x_{1}^* , \ldots, x_{n-1}^*, x_{n}^* \Big)
&\mbox{for $N=2n+1$},
\\[6pt]
{1 \over r} ~\Big( q^{-\rho_{n}} ~x_{n} , q^{-\rho_{n-1}} ~x_{n-1},
\ldots, q^{-\rho_{1}} ~x_{1} , x_{1}^* , \ldots, x_{n-1}^*, x_{n}^* \Big)
&\mbox{for $N=2n$},
\end{array}\right.
\end{equation}
which is clearly normalized to $1$ : $\hs{\perp\!}{\!\perp} = r^2 /
r^2= 1$, from relation (\ref{radius}). Thus
the matrix valued element $e_{\perp} \in \Mat_{N}(S^{N-1}_q)$ given by
\begin{equation}                                             \label{norbun}
e_{\perp} = \ket{\!\perp} \bra{\perp\!}
\end{equation}
is a self-adjoint idempotent, i.e.  $e_{\perp}^2 = e_{\perp}$,
$e_{\perp}^* = e_{\perp}$.  The corresponding finite
projective module over the sphere $S^{N-1}_q$ will be named the {\it
normal bundle} over the sphere $S^{N-1}_q$.  Furthermore, the
self-adjoint idempotent
\begin{equation}                                         \label{cotbun}
e_{Cot} = \b{I}_N - \ket{\!\perp} \bra{\perp\!}
\end{equation}
is a natural candidate for the {\it cotangent bundle} over the sphere
$S^{N-1}_q$. \\
Again, the rank of the
bundles is computed by
combining the matrix trace with the $1$-dimensional representation in
(\ref{1reps}),
${\rm rank}(e) = \tau_\lambda \circ Tr(e)$.  
Then, straightforward computations give that,
as expected, the projector $e_{\perp}$ is of rank $1$ while the projector 
$e_{Cot}$ is of rank $N-1$,
\begin{eqnarray}
&& {\rm rank}(e_{\perp}) = \tau_\lambda \circ Tr(e_{\perp}) 
= {1\over r^2} ~ \tau_\lambda \Big( q^{-2\rho_n} x_n^* x_n + \cdots + x_n x_n^* \Big) =
1 ,
\nonumber \\ && {\rm rank}(e_{Cot}) = \tau_\lambda \circ Tr(e_{Cot}) = N-1 .
\end{eqnarray}
Relations of the cotangent projector (\ref{cotbun}) with differential calculi on the spheres
$S^{N-1}_q$ will be analyzed elsewhere.

\initiate
\section{Final remarks}

We have presented a description of (the algebra of functions on) the
quantum Euclidean spheres
$S^{N-1}_q$ by means of projections $e\in\Mat_{2^n}(S^{N-1}_q)$, with 
$N=2n$ or $N=2n+1$
which, in turn, can be regarded as bundles over the spheres $S^{N-1}_q$.
Apart from some identifications in `lower dimensions', 
these spheres are different from analogous objects described in
\cite{Wor87,Pod87,VakSoi90,ConLan00,DabLanMas00,BonCicTar00,Sit01,BrzGon01,DabLan01,
HonSzy01,Var01,ConDub-Vio01,AscBon01}. The monopole presented here on the sphere
$S^2_q$  coincides with the one constructed in \cite{BrzMaj00}
on the equator spheres of \cite{Pod87}. Our $q$-solitons and 
$q$-instantons on the
spheres $S^3_q$, $S^4_q$ seems to be different from
analogous objects recently found in
\cite{DabLanMas00,BonCicTar00,Sit01,BrzGon01,DabLan01}. It is also clearly
different from the instanton constructed in \cite{ConLan00}.
The present Euclidean spheres  are characterized by a `homological
dimension drop' which signals the fact that, contrary to the spheres 
$S^4_{\theta}$ of
\cite{ConLan00}, they are {\em not} noncommutative manifolds and their geometry
cannot be a solution of  homological equations like the ones in \cite{ConLan00}.
In \cite{ConLan00} the spheres $S^4_{\theta}$ were endowed with a 
noncommutative
geometry via an even spectral triple $(\c{A},\c{H},D,\gamma)$ where $\c{A}$
is a noncommutative algebra with involution * acting on a 
$\IZ_2$-graded Hilbert space
$\c{H}$ (with grading given by $\gamma$) while $D$ is a self-adjoint 
operator on $\c{H}$
with compact resolvent and such that $[D,a]$ is bounded for any 
$a\in\c{A}$ \cite{Con94}.
The operator D specify both the metric on the state space of $\c{A}$ 
and the K-homology
fundamental class \cite{Con96}. The geometry for the spheres $S^4_{\theta}$ is 
constructed by deforming
the commutative triple $(C^{\infty}(S^4),\c{H},D,\gamma_5)$, where 
$D$ is the Dirac
operator on the Hilbert space $\c{H}$ of square integrable spinors 
over $S^4$ and the
grading is given by the `fifth gamma matrix'. In fact, one has an isospectral
deformation since both $\c{H}$ and $D$ are kept fixed, so that all 
spectral data of the
geometry are unchanged, while the algebra and its representations are 
deformed. 

An important problem that we leave for the future is the computation of topological
charges, notably on the quantum Euclidean sphere $S^4_q$ and higher dimensional even
spheres. This is a difficult task since it involves the construction of
Fredholm modules and their characters  which, via the noncommutative index theorem of
\cite{Con94}, pair integrally with the $K_0$ group. That the construction is
difficult it is  already evident from the analogous constructions
\cite{MasNakWat90,MasNakWat91} for the spheres
$S^2_q$ and $S^3_q$. 
As a preliminary step for this construction we introduce some
representations of the algebra of the quantum Euclidean sphere $S^4_q$ (similar
representations can be constructed for any sphere and generalize  known results for
the quantum  spheres $S^2_q$ and $S^3_q$). 

Let us then consider again the algebra
of $S^4_q$ which is specified by the commutation relations (\ref{cr4}).  
And  let $\c{H}$ be an infinite dimensional Hilbert space with
orthonormal basis $\{\psi_{n,m} , ~n,m = 0,1, 2,
\cdots ~\}$. We have two representations (see also
\cite{Fio95})
\begin{eqnarray}\label{reps4plus}
&& \pi(x_0) \psi_{n,m} = \pm ~r ~\sqrt{ { 1 + q \over 1+q^3}}
~q^{n + 1/2} ~q^{m + 1} ~\psi_{n,m},
\nn \\ && \nn \\
&& \pi(x_1) \psi_{n,m} = ~r ~\sqrt{ { 1 -
q^{2(n+1)}  \over
1+q^3}} ~q^{m+1} ~\psi_{n+1,m},
\nn \\
&& \nn \\
&& \pi(x_1^*) \psi_{n,m} = ~r ~\sqrt{ {
1 - q^{2n}  \over
1+q^3}} ~q^{m+1}
~\psi_{n-1,m} ,
\nn \\
&& \nn \\
&& \pi(x_2) \psi_{n,m} = ~r ~\sqrt{ { 1 - q^{2(m+1)}
\over 1+q^3}} ~ \psi_{n,m+1} ,
\nn \\
&& \nn \\
&& \pi(x_2^*)
\psi_{n,m} = ~r ~\sqrt{ { 1 - q^{2m}  \over 1+q^3}} 
~\psi_{n,m-1} .
\end{eqnarray}
We notice that for $q<1$ any power of the operator $\pi(x_0)$
is a trace class operator while this is not the
case for the operators $\pi(x_1)$, $\pi(x_1^*)$,
$\pi(x_2)$ and $\pi(x_2^*)$ nor
for any of their powers. Were this not the case there would be a
contradiction with the algebra relations (\ref{cr4}). A detailed
analysis of the previous representations will be reported somewhere else.

\initiate
\section{Acknowledgment}

G.L. thanks Genevi{\`e}ve for her wonderful hospitality at the M-PSC in Orsay
where this work initiated. We are grateful to Alain Connes, Tomasz 
Brzezinski, Ludwik Dabrowski, 
Gaetano Fiore, Harold Steinacker and Roberto Trinchero for discussions and
suggestions.

\vfill\eject


\begin{thebibliography}{10}

\bibitem{AscBon01} P.~Aschieri, F.~Bonechi,
``On the Noncommutative Geometry of Twisted Spheres,''
{\em Lett. Math. Phys.} {\bf 59} (2002) 133,
\href{http://xxx.lanl.gov/abs/math.math.QA/0108136}{{\tt
math.QA/0108136}}.

\bibitem{BonCicTar00}
F.~Bonechi, N.~Ciccoli, and M.~Tarlini, ``Noncommutative Instantons on the
     4-sphere from quantum groups,''
     {\em Commun. Math. Phys.} {\bf 226} (2002) 419,
\href{http://xxx.lanl.gov/abs/math.QA/0012236}{{\tt math.QA/0012236}}.

\bibitem{BrzGon01}
T.~Brzezinski and C.~Gonera, ``Non-commutative 4-spheres based on all Podle\'s
     2-spheres and beyond,'' 
{\em Lett. Math. Phys.} {\bf 54} (2000) 315,
\href{http://xxx.lanl.gov/abs/math.QA/0101129}{{\tt math.QA/0101129}}.

\bibitem{BrzMaj98} T.~Brzezinski and S.~Majid, 
``Line bundles on quantum spheres,'' In: Particles, Fields and
Gravitation. J. Rembielinski (ed), AIP Woodbury, New York, pp.3-8 (1998),
\href{http://xxx.lanl.gov/abs/math.QA/9803003}{{\tt math.QA/9807052}}.

\bibitem{BrzMaj00} T.~Brzezinski and S.~Majid,
``Quantum geometry of algebra factorisations and coalgebra bundles,''
{\em Commun. Math. Phys.} {\bf 213} (2000) 491,
\href{http://xxx.lanl.gov/abs/math.QA/9803003}{{\tt math.QA/9808067}}.
 
\bibitem{CerFioMad00a}
B.L.~Cerchiai, G.~Fiore, and J.~Madore, ``Geometrical tools for quantum
     euclidean spaces,'' {\em Commun.\ Math.\ Phys.} {\bf 217} (2001)
521, \href{http://xxx.lanl.gov/abs/math.QA/0002007}{{\tt math.QA/0002007}}.
 
\bibitem{Con94} A. Connes, {\it Noncommutative geometry}. Academic Press 1994.

\bibitem{Con96} A. Connes, ``Gravity coupled with matter and
foundation of  noncommutative geometry,'' {\em Commun. Math. Phys.}
{\bf 182} (1996) 155.

\bibitem{ConDub-Vio01} A.~Connes, M.~Dubois-Violette,
``Noncommutative finite-dimensional manifolds. I. Spherical manifolds and related 
 examples,''
\href{http://xxx.lanl.gov/abs/math.math.QA/0107070}{{\tt math.QA/0107070}}.

\bibitem{ConLan00}
A.~Connes and G.~Landi, ``Noncommutative manifolds the instanton algebra and
     isospectral deformations,'' {\em Commun.\ Math.\ Phys.}, 
{\bf 221} (2001) 141,
     \href{http://xxx.lanl.gov/abs/math.QA/0011194}{{\tt math.QA/0011194}}.

\bibitem{DabLan01}
L.~Dabrowski and G.~Landi, ``Instanton algebras and quantum 4-spheres,''
{\em  Differ. Geom. Appl., in press},
\href{http://xxx.lanl.gov/abs/math.QA/0101177}{{\tt math.QA/0101177}}.
 
\bibitem{DabLanMas00}
L.~Dabrowski, G.~Landi, and T.~Masuda, ``Instantons on the quantum 4-spheres
     {$S_q^4$},'' {\em Commun.\ Math.\ Phys.} {\bf 221} (2001) 161,
\href{http://xxx.lanl.gov/abs/math.QA/0012103}{{\tt math.QA/0012103}}.
 
\bibitem{FadResTak89}
L.~D. Faddeev, N.~Y. Reshetikhin, and L.~A. Takhtajan, ``Quantization of {L}ie
     groups and {L}ie algebras,'' {\em Lenin. Math. Jour.} {\bf 1} (1990) 193.

\bibitem{Fio95}
G.~Fiore,  ``The Euclidean Hopf algebra $U_q(e^N)$ and its fundamental
Hilbert space representations,''  {\em J. Math. Phys.} {\bf 36} (1995) 4363,
\href{http://xxx.lanl.gov/abs/hep-th/9407195}{{\tt hep-th/9407195}}. \\
G.~Fiore, ``The q-Euclidean algebra $U_q(e^N)$ and the corresponding
q-Euclidean lattice,'' {\em Int. J. Mod. Phys.} {\bf A11} (1996) 863,
\href{http://xxx.lanl.gov/abs/q-alg/9506028}{{\tt q-alg/9506028}}.

\bibitem{Gra00} J.M.~Gracia-Bondia, J.C.~Varilly, H.~Figueroa,
{\it Elements of noncommutative geometry}. Birkhauser, 2000.

\bibitem{HajMaj98}
P.~Hajac and S.~Majid, ``Projective module description of the $q$-monopole,''
     {\em Commun. Math. Phys.} {\bf 206} (1999) 247,
\href{http://xxx.lanl.gov/abs/math.QA/9803003}{{\tt math.QA/9803003}}.

\bibitem{HonSzy01} J.H.~Hong and W.~Szyma{\'n}ski, ``Quantum spheres and
projective spaces as graph algebras,'' University of Newcastle Preprint, January
2001.

\bibitem{Lan97} G.~Landi, {\it An introduction to noncommutative
spaces and their
geometries}. Springer~1997.

\bibitem{Lan00} G. Landi, ``Deconstructing Monopoles and Instantons,''
{\em Rev. Math. Phys.} {\bf 12} (2000) 1367,
\href{http://xxx.lanl.gov/abs/math-ph/9812004}{{\tt math-ph/9812004}}.

\bibitem{Mad99}
J.~Madore, {\em An introduction to noncommutative differential geometry and its
     physical applications}.
\newblock No.~257 in London Mathematical Society Lecture Note Series. Cambridge
     University Press, second~ed. , 1999.

\bibitem{MasNakWat90}
T.~Masuda, Y.~Nakagami and J.~Watanabe, ``Noncommutative differential 
geometry on the quantum $SU(2)$. I: An algebraic 
viewpoint,'' {\em $K$-Theory} {\bf 5} (1990) 157.

\bibitem{MasNakWat91}
T.~Masuda, Y.~Nakagami and J.~Watanabe, ``Noncommutative differential 
geometry on the quantum two sphere of P.~Podle{\'s}. I: An algebraic 
viewpoint,'' {\em $K$-Theory} {\bf 5} (1991) 151.

\bibitem{Ogi92a}
O.~Ogievetsky, ``Differential operators on quantum spaces for {$GL_q(n)$} and
     {$SO_q(n)$},'' {\em Lett.\ Math.\ Phys.} {\bf 24} (1992) 245.

\bibitem{Pod87}
P.~Podle{\'s}, ``Quantum spheres,'' {\em Lett.\ Math.\ Phys.} {\bf 14} (1987)
     193.

\bibitem{Sit01}
A.~Sitarz, ``More noncommutative 4-spheres,''
{\em Lett.  Math.  Phys.} {\bf 55} (2001) 127,      
\href{http://xxx.lanl.gov/abs/math-ph/0101001}{{\tt math-ph/0101001}}.

\bibitem{VakSoi90} L.L.~Vaksman and Y.S.~Soibelman, ``Algebra of functions on
quantum $SU(n+1)$ group and odd dimensional quantum spheres,''
{\em Algebra-i-Analiz} {\bf 2} (1990) 101.

 \bibitem{Var01} J.C.~Varilly, ``Quantum symmetry groups of 
noncommutative spheres,'' {\em Commun. Math. Phys.} {\bf 221} (2001) 
511,
\href{http://xxx.lanl.gov/abs/math.math.QA/0102065}{{\tt math.QA/0102065}}.

\bibitem{Wor87} S.L.~Woronowicz, {\it Twisted $SU(2)$ group. An 
example of a noncommutative differential calculus}. Publ. Res. Inst. Math. 
Sci. 23 (1987) 117.
 
\end{thebibliography}

\providecommand{\href}[2]{#2}\begingroup\raggedright\endgroup

\end{document}